\documentclass[10pt]{article}

\usepackage{amsmath,amsfonts, latexsym,graphicx}

\newcommand{\U}{{\mathcal U}}
\newcommand{\V}{{\mathcal V}}
\newcommand{\0}{{\mathbf 0}}
\newcommand{\C}{{\mathbb C}}

\newcommand{\R}{{\mathbb R}}

\newcommand{\N}{{\mathbb N}}

\newcommand{\D}{{\mathbb D}}
\newcommand{\W}{{\mathcal W}}
\newcommand{\strat}{{\mathfrak S}}

\newcommand{\ob}{{{}^{{}^\circ}\hskip-.04in B}}
\newcommand{\od}{{{}^{{}^\circ}\hskip-.02in \D}}

\newcommand{\tr}{{\operatorname{tr}}}

\newtheorem{defn0}{Definition}[section]
\newtheorem{prop0}[defn0]{Proposition}
\newtheorem{conj0}[defn0]{Conjecture}
\newtheorem{thm0}[defn0]{Theorem}
\newtheorem{lem0}[defn0]{Lemma}
\newtheorem{corollary0}[defn0]{Corollary}
\newtheorem{example0}[defn0]{Example}
\newtheorem{remark0}[defn0]{Remark}
\newtheorem{question0}[defn0]{Question}

\newenvironment{defn}{\begin{defn0}\hskip -.06in .}{\end{defn0}}
\newenvironment{prop}{\begin{prop0}\hskip -.06in .}{\end{prop0}}

\newenvironment{thm}{\begin{thm0}\hskip -.06in .}{\end{thm0}}
\newenvironment{lem}{\begin{lem0}\hskip -.06in .}{\end{lem0}}
\newenvironment{cor}{\begin{corollary0}\hskip -.06in .}{\end{corollary0}}
\newenvironment{exm}{\begin{example0}\hskip -.06in .\rm}{\end{example0}}
\newenvironment{rem}{\begin{remark0}\hskip -.06in .\rm}{\end{remark0}}
\newenvironment{ques}{\begin{question0}\hskip -.06in .\rm}{\end{question0}}

\newcommand{\defref}[1]{Definition~\ref{#1}}
\newcommand{\propref}[1]{Proposition~\ref{#1}}
\newcommand{\thmref}[1]{Theorem~\ref{#1}}
\newcommand{\lemref}[1]{Lemma~\ref{#1}}
\newcommand{\corref}[1]{Corollary~\ref{#1}}
\newcommand{\exref}[1]{Example~\ref{#1}}
\newcommand{\secref}[1]{Section~\ref{#1}}
\newcommand{\remref}[1]{Remark~\ref{#1}}

\newcommand{\qed}{\mbox{$\Box$}}
\newenvironment{proof}{\noindent {\bf Proof.}}{\qed\vskip 6pt}

\title{Real Analytic Milnor Fibrations and a Strong  \L ojasiewicz Inequality\footnote{\mbox{   }   AMS subject classifications 32C18, 14P15, 26E05, 32B10, 32B20.
\newline   \mbox{   } \mbox{   }   Keywords:  \L ojasiewicz inequality, real Milnor fibration, Thom stratifications}}

\author{David B. Massey}

\date{}

\begin{document}

\baselineskip = 14pt

\maketitle

\begin{abstract} We give a strong version of a classic inequality of \L ojasiewicz; one which collapses to the usual inequality in the complex analytic case. We show that this inequality for real analytic functions allows us to construct a real Milnor fibration inside a ball. \end{abstract}

\sloppy

\section{Introduction}\label{sec:intro}

Suppose that $\U$ is an open neighborhood of the origin in $\C^N$, and that $f_{{}_\C}:(\U, \0)\rightarrow (\C, 0)$ is a complex analytic function.

In the now-classic book \cite{milnorsing}, Milnor shows that one has what is now called the Milnor fibration of $f_{{}_\C}$ (at $\0$). The Milnor fibration is THE fundamental device in the study of the topology of the hypersurface $X$ defined by the vanishing of $f_{{}_\C}$.

In fact, there are two Milnor fibrations associated with $f_{{}_\C}$: one defined on small spheres, and one defined inside small open balls. Both of these are referred to as {\bf the} Milnor fibration because the two fibrations are diffeomorphic. We wish to be precise. 

For $\epsilon>0$, let $S_\epsilon$ (resp., $B_\epsilon$, $\ob_\epsilon$) denote the sphere (resp., closed ball, open ball) of radius $\epsilon$, centered at the origin. In the special case of balls in $\R^2\cong\C$, we write $\D_\epsilon$ in place of $B_\epsilon$ and $\od_\epsilon$ in place of $\ob_\epsilon$. Finally, we write $\D_\epsilon^*$ in place of $\D_\epsilon-\{\0\}$, and $\od_\epsilon^*$ in place of $\od_\epsilon-\{\0\}$.

One version of the Milnor fibration (at the origin) is given by: there exists $\epsilon_0>0$ such that, for all $\epsilon$ such that $0<\epsilon\leq\epsilon_0$, the map $f_{{}_\C}/|f_{{}_\C}|:S_\epsilon-S_\epsilon\cap X\rightarrow S^1\subseteq \C$ is a smooth, locally trivial ly fibration, whose diffeomorphism-type is independent of the choice of $\epsilon$ (see \cite{milnorsing}). 

The second, diffeomorphic version of the Milnor fibration is given by: there exists $\epsilon_0>0$ such that, for all $\epsilon$ such that $0<\epsilon\leq\epsilon_0$, there exists a $\delta_0>0$ such that, for all $\delta$ such that $0<\delta\leq\delta_0$, the map $f_{{}_\C}:\ob_\epsilon\cap {f_{{}_\C}}^{-1}(\partial\D_\delta)\rightarrow \partial\D_\delta\cong S^1\subseteq \C$ is a smooth, locally trivial ly fibration, whose diffeomorphism-type is independent of the choice of $\epsilon$ and (sufficiently small) $\delta$ (see Theorem 5.11 of \cite{milnorsing} and  \cite{relmono}). The primary advantage to this second characterization of the Milnor fibration is that it compactifies nicely to yield a locally trivial  fibration $f_{{}_\C}:B_\epsilon\cap {f_{{}_\C}}^{-1}(\partial\D_\delta)\rightarrow \partial\D_\delta$, which, up to homotopy, is equivalent to either of the two previously-defined Milnor fibrations.

We shall not summarize the important properties of the Milnor fibration here, but refer the reader to \cite{milnorsing}, \cite{dimcasing}, \cite{randell}, and the Introduction to \cite{lecycles}. We wish to emphasize that our discussion of the Milnor fibration above assumes that $f_{{}_\C}$ is a complex analytic function.

\bigskip

Of course, a complex analytic function yields a pair of real analytic functions, coming from the real and imaginary parts of the complex function, and one can ask the more general question: {\bf when does a pair of real analytic functions possess one or both types of Milnor fibrations?}

\bigskip

This topic of real analytic Milnor fibrations is complicated and interesting, and gives rise to many questions.

\smallskip

In Chapter 11 of \cite{milnorsing}, Milnor discusses, fairly briefly, some results in the case of the real analytic function $f=(g, h)$. He considers the very special case where $f$ has an isolated critical point at the origin, and shows that, while the restriction of $f$ still yields a fibration over a small circle inside the ball, $f/|f|$ does not necessarily yield the projection map of a fibration from $S_\epsilon-S_\epsilon\cap X$ to $S^1$; see p. 99 of \cite{milnorsing}. Can one relax the condition that $f$ has an isolated critical point and still obtain a locally trivial  fibration $f:\ob_\epsilon\cap f^{-1}(\partial\D_\delta)\rightarrow \partial\D_\delta$? Are there reasonable conditions that guarantee that $f/|f|:S_\epsilon-S_\epsilon\cap X\rightarrow S^1$ is a locally trivial  fibration which is diffeomorphic to the fibration inside the ball? Such questions have been investigated by a number of researchers; see \cite{milnorsing}, \cite{jacquemard}, \cite{pichonseade}, \cite{seadebook}, \cite{ruassantos}, and \cite{santos}.

\medskip

An obvious approach to answering the question about the existence of real analytic Milnor fibrations is simply to try to isolate what properties of a complex analytic function are used in proofs that MIlnor fibrations exist. We write ``proofs'' and not ``the proof'' here because we will not follow Milnor's proof, but rather follow L\^e's proof in \cite{relmono} of the existence of Milnor fibrations inside the ball for complex analytic functions.

L\^e's proof consists almost solely of using the existence of a Thom (or $a_f$, or {\it good}) stratification. In \cite{hammlezariski}, Hamm and L\^e, following a suggestion of Pham, used the complex analytic  \L ojasiewicz inequality (see \corref{cor:complexloj} below) and a ``trick'' to show that Thom stratifications exist.

\bigskip

Our goal in this paper is very modest: we will give the ``correct'' generalization of the complex analytic  \L ojasiewicz inequality, and then show that if a pair (or quadruple, or octuple) of real analytic functions satisfies this new {\it strong  \L ojasiewicz inequality}, then the Milnor fibration inside a ball exists.

\bigskip

In the remainder of the introduction, we will summarize our primary definition and result.

\bigskip

 Let $\U$ now denote an open subset of  $\R^n$, and $p$ denote a point in $\U$. Let $g$ and $h$ be real analytic functions from $\U$ to $\R$, and let $f:=(g, h):\U\rightarrow \R^2$. Recall the classic inequality of  \L ojasiewicz \cite{lojbooknotes} (see also \cite{biermilsemi}).

\smallskip

\begin{thm} \label{thm:loj} There exists an open neighborhood $\W$ of $p$ in $\U$, and $c,\theta\in\R$ such that $c>0$, $0<\theta<1$, and, for all $x\in\W$,
$$
|g(x)-g(p)|^\theta\leq c\big |\nabla g(x)\big|,
$$
where $\nabla g$ is the gradient vector.
\end{thm}

\begin{rem}\label{rem:noc} The phrasing above is classical, and convenient in some arguments. However, one can also fix the value of the constant  $c$ above to be any $c>0$, e.g., $c=1$. In other words, one may remove the reference to $c$ in the statement \thmref{thm:loj} and simply use the inequality
$$
|g(x)-g(p)|^\theta\leq \big |\nabla g(x)\big|.
$$
The argument is easy, and simply requires one to pick a larger $\theta$ (still less than $1$). As we shall not use this ``improved'' statement, we leave the proof as an exercise.
\end{rem}

\medskip

\thmref{thm:loj} implies a well-known complex analytic version of itself. One can easily obtain this complex version by replacing $g$ by the square of the norm of the complex analytic function. However, we shall prove two not so well-known more general corollaries, \corref{cor:lojtuple} and \corref{cor:lojtuple2}, which yield the complex analytic statement. 

\medskip

If $n$ is even, say $n=2m$, then we may consider the complexified version of $f$ by defining $f_{{}_\C}$ by
$$
f_{{}_\C}(x_1+iy_1, \dots, x_m+iy_m):= g(x_1, y_1, \dots, x_m, y_m) + ih(x_1, y_1, \dots, x_m, y_m).
$$

From \corref{cor:lojtuple2}, we immediately obtain:

\begin{cor}\label{cor:complexloj} Suppose that $f_{{}_\C}$ is complex analytic. Then, there exists an open neighborhood $\W$ of $p$ in $\U\subseteq\C^n$, and $c,\theta\in\R$ such that $c>0$, $0<\theta<1$, and, for all $z\in\W$,
$$
|f_{{}_\C}(z)-f_{{}_\C}(p)|^\theta\leq c\big |\nabla f_{{}_\C}(z)\big|,
$$
where $\nabla f_{{}_\C}$ denotes the complex gradient.
\end{cor}

\bigskip

Our generalization of this complex analytic \L ojasiewicz inequality is:

\medskip

\begin{defn} We say that $f=(g,h)$ satisfies the {\bf strong \L ojasiewicz inequality at $p$} or that {\bf $f$ is \L-analytic at $p$} if and only if 
there exists an open neighborhood $\W$ of $p$ in $\U$, and $c,\theta\in\R$ such that $c>0$, $0<\theta<1$, and, for all $x\in\W$,
$$
|f(x)-f(p)|^\theta\leq c\cdot\underset{|(a, b)|=1}{\operatorname{min}}\big |a\nabla g(x)+b\nabla h(x)\big|.
$$
\end{defn}

\medskip

\noindent {\bf Main Result}. {\it Suppose that  $f(\0)=\0$, that $f$ is not locally constant near the origin, and that $f$ is \L-analytic at $\0$.

Then, for all $0<\delta\ll\epsilon\ll 1$, $f:B_\epsilon\cap f^{-1}(\partial\D_\delta)\rightarrow \partial\D_\delta$ is a proper, stratified submersion, and so $f:B_\epsilon\cap f^{-1}(\partial\D_\delta)\rightarrow \partial\D_\delta$ and $f:\ob_\epsilon\cap f^{-1}(\partial\D_\delta)\rightarrow \partial\D_\delta$ are locally trivial  fibrations.

Moreover, the diffeomorphism-type of $f:\ob_\epsilon\cap f^{-1}(\partial\D_\delta)\rightarrow \partial\D_\delta$ is independent of the appropriately small choices of $\epsilon$ and $\delta$.
}

\bigskip

We have presented the main definition and result above in the case of real analytic functions into $\R^2$; this was for simplicity of the discussion. In fact,  in \defref{def:strongl} and \corref{cor:main}, we give our main definition and result when the dimension of the codomain is arbitrary.

\smallskip

We thank T. Gaffney for pointing out the existence of nap-maps (see \defref{def:napmap}) in all dimensions.

\section{Singular Values of Matrices}\label{sec:singval}

In this section, we wish to recall some well-known linear algebra, and establish/recall some inequalities that we will need in later sections.

Let $\vec v_1, \dots \vec v_k$ be vectors in $\R^n$. Let $A$ denote the $n\times k$ matrix whose $i$-th column is $\vec v_i$. Let $M$ denote the $k\times k$ matrix $A^tA$.
Consider the function $n$ from the unit sphere centered at the origin in $\R^k$ into the non-negative real numbers given by $n(t_1, \dots, t_k) := |t_1\vec v_1 + \dots +t_k \vec v_k|$. The critical values of $n$ are the {\it singular values of $A$}, and they are the square roots of the necessarily non-negative eigenvalues of $M$. It is traditional to index the singular values in a decreasing manner, i.e., we let the singular values of $A$  (which need not be distinct) be denoted by $\sigma_1, \dots, \sigma_k$, where $\sigma_1\geq\dots\geq\sigma_k$. We denote the eigenvalues of $M$ by $\lambda_i:=\sigma_i^2$.  The singular value $\sigma_1$ is the norm of $A$. The minimum singular value $\sigma_k$ will be of particular interest throughout this paper. Note that
$$
\sigma_k= \underset{|(t_1, \dots, t_k)|=1}{\operatorname{min}}\big |t_1\vec v_1+\dots +t_k\vec v_k\big|.\leqno{(\dagger)}
$$
 Also note that the trace of $M$, $\tr M$, is equal to $|\vec v_1|^2+\dots+|\vec v_k|^2$, which is equal to $\lambda_1+\dots+\lambda_k=\sigma_1^2+\dots+\sigma_k^2$. Hence, $\tr M=0$ if and only if $A=0$ if and only if $\sigma_1=0$. It will be important for us later that the singular values of $A$ vary continuously with the entries of $A$. Also, as the non-zero eigenvalues of $A^tA$ and $AA^t$ are the same, the non-zero singular values of $A$ are equal to the non-zero singular values of $A^t$; in particular, $\sigma_1(A)=\sigma_1(A^t)$, i.e., the norm of a matrix is equal to the norm of its transpose.
 
 \smallskip

The following proposition contains the fundamental results on singular values that we shall need; these results are known, but we include brief proofs for the convenience of the reader.

\begin{prop}\label{prop:singval}
\begin{enumerate}
\item $\sigma_k=0$ if and only if $\vec v_1, \dots, \vec v_k$ are linearly dependent;
\item if $|\vec v_1|^2+\dots+|\vec v_k|^2\neq 0$, then 
$$\frac{1}{\sqrt{k}}\leq\frac{\sigma_1}{\sqrt{|\vec v_1|^2+\dots+|\vec v_k|^2}}\leq 1;$$
and
\item $k(\det M)^{1/k}\leq |\vec v_1|^2+\dots+|\vec v_k|^2$, with equality holding if and only if $\vec v_1, \dots, \vec v_k$ all have the same length and are pairwise orthogonal.
\end{enumerate}
\end{prop}
\begin{proof} Item 1 is immediate from $(\dagger)$. 

Assume that $ |\vec v_1|^2+\dots+|\vec v_k|^2\neq 0$. 
As $|\vec v_1|^2+\dots+|\vec v_k|^2 = \sigma_1^2+\dots+\sigma_k^2$, the right-hand inequality in Item 2 is immediate. For all $i$, let 
$$
r_i:=\frac{\sigma_i}{\sqrt{|\vec v_1|^2+\dots+|\vec v_k|^2}}.
$$
Then, $r_1^2+\dots+r_k^2 =1$ and $r_1^2\geq\dots\geq r_k^2$. It follows that $r_1^2\geq 1/k$; this proves the left-hand inequality in Item 2.

The inequality in Item 3 is nothing more than the fact that the geometric mean of the eigenvalues of $M$ is less than or equal to the arithmetic mean, where we have again used that the trace of $M$ is $|\vec v_1|^2+\dots+|\vec v_k|^2$. In addition, these two means are equal if and only if all of the eigenvalues of $M$ are the same. This occurs if and only if all of singular values of $A$ are the same, which would mean that the function  $n(t_1, \dots, t_k) := |t_1\vec v_1 + \dots +t_k \vec v_k|$ is constant on the unit sphere. The remainder of Item 3 follows easily.
\end{proof}

\medskip

We need some results for dealing with compositions. Let $\sigma_i(C)$  denote the $i$-th singular value of a matrix $C$, indexed in descending order. Let $\lambda_i(C):=\sigma_i^2(C)$, i.e., $\lambda_i(C)$ is the $i$-th eigenvalue of $C^tC$, where the eigenvalues are indexed in descending order. 

\smallskip Let $B$ be an $m\times n$ matrix.

\begin{prop}\label{prop:prodsv} 
$$
\sigma_k(BA)\geq \sigma_n(B)\sigma_k(A).
$$
Hence,
$$
\left[\det(A^tB^tBA)\right]^{1/k}\geq \sigma_k^2(BA)\geq \sigma_n^2(B)\sigma_k^2(A).
$$
\end{prop}
\begin{proof} Let $\vec u$ be a unit vector in $\R^k$, written as a $k\times 1$ matrix. If $A\vec u=0$, then $\sigma_k(A)=0$ and the first inequality is immediate. So, assume that $A\vec u\neq 0$. Then,
$$
|BA\vec u|=\left| B\left(\frac{A\vec u}{|A\vec u|}\right)\right|\cdot  |A\vec u|\cdot  \geq \sigma_n(B)\sigma_k(A).
$$
This proves the first inequality.

As $\left[\det(A^tB^tBA)\right]^{1/k}$ is the geometric mean of the $\sigma_i^2(BA)$, the second set of inequalities are immediate.
\end{proof}

\begin{lem}\label{lem:traceprod} Suppose that $P$ and $Q$ are $r\times r$ real matrices, and that $P$ is diagonalizable (over the reals) and has no negative eigenvalues.  Then,
$$
\lambda_r(P)\tr(Q)\leq \tr(PQ)\leq \lambda_1(P)\tr(Q).
$$
\end{lem}
\begin{proof} Suppose that $P=S^{-1}DS$, where $D$ is diagonal. Then, $\tr(PQ)=\tr(S^{-1}DSQ)=\tr(DSQS^{-1})$. Let $B=SQS^{-1}$. Hence, $\tr(B)=\tr(Q)$, and 
$$
\tr(PQ)=\tr(DB)=\sum_{i=1}^r(DB)_{i,i}=\sum_{i=1}^r\sum_{j=1}^r D_{i, j}B_{j, i}= \sum_{i=1}^r D_{i,i}B_{i,i}.
$$
As $\lambda_r(P)\leq D_{i,i}\leq\lambda_1(P)$, the conclusion follows.
\end{proof}

\begin{prop}\label{prop:traceprod} 
$$
n\sigma^2_n(A^t)\sigma^2_n(B)\leq \tr(A^tB^tBA)\leq n\sigma^2_1(A)\sigma^2_1(B).
$$
\end{prop}
\begin{proof} As $\tr(A^tB^tBA)=\tr(AA^tB^tB)$, we may apply \lemref{lem:traceprod} with $P=A^tA$ and $Q=B^tB$, and conclude that
$$
\sigma^2_n(A^t)[n\sigma^2_n(B)]\leq\sigma^2_n(A^t)\tr(B^tB)\leq \tr(A^tB^tBA)\leq \sigma^2_1(A^t)\tr(B^tB)\leq \sigma^2_1(A^t)[n\sigma^2_1(B)].
$$
Now, use that $\sigma_1(A)=\sigma_1(A^t)$.
\end{proof}

\begin{cor}\label{cor:slw} If $BA\neq 0$, then $\tr(A^tB^tBA)$, $\sigma_1^2(B)$, and $\sigma_1^2(A)$ are not zero and 
$$
\frac{k\left[\det(A^tB^tBA)\right]^{1/k}}{\tr(A^tB^tBA)}\geq\frac{k\sigma_n^2(B)\sigma_k^2(A)}{n\sigma_1^2(B)\sigma_1^2(A)}.
$$
\end{cor}
\begin{proof} Combine \propref{prop:traceprod} with \propref{prop:prodsv}.
\end{proof}

\section{\L-maps and \L-weights}\label{sec:lanal}

As in the introduction, let $\U$ be an open subset of $\R^n$. Let $f=(g,h)$ be a $C^1$ map from $\U$ to $\R^2$ and let $G:=(g_1, \dots, g_k)$ be a $C^1$ map from $\U$ to $\R^k$. By the {\it critical locus of $G$}, $\Sigma G$, we mean the set of points where $G$ is not a submersion (this is reasonable since our main hypothesis in most results will imply that $n\geq k$ or that $G$ is locally constant). A number of our results will apply only with the stronger assumption that $f$ and $G$ are real analytic, but we shall state that hypothesis explicitly as needed.

Throughout this section, we let $A=A(x)$ denote the $n\times k$ matrix which has the gradient vector of $g_i(x)$ as its $i$-th column (i.e., $A$ is the transpose of the derivative matrix $[d_xG]$ of $(g_1, \dots, g_k)$), and let $M:=A^tA$. We will be applying the results of \secref{sec:singval} to $A$ and $M$. Let $\sigma_1(x), \dots, \sigma_k(x)$ denote the singular values of $A(x)$, indexed in decreasing order.

\medskip

We have the following corollary to \thmref{thm:loj}

\begin{cor} \label{cor:lojtuple} Suppose that $G$ is real analytic. Then, there exists an open neighborhood $\W$ of $p$ in $\U$, and $c,\theta\in\R$ such that $c>0$, $0<\theta<1$, and, for all $x\in\W$,
$$
|G(x)-G(p)|^\theta\leq c \sqrt{|\nabla g_1(x)|^2+\dots+|\nabla g_k(x)|^2}.
$$
\end{cor}
\begin{proof} For notational convenience, we shall prove the result for $k=2$, using $f=(g, h)$, and we shall assume that $f(p)=0$; the proof of the general case proceeds in exactly the same manner. 

\smallskip

We will prove that there exists $\W$, $c$, and $\theta$ as in the statement such that the inequality holds at all $x\in\W$ such that $f(x)\neq 0$. This clearly suffices to prove the corollary. We will also assume that $|\nabla g(p)|^2+|\nabla h(p)|^2= 0$ for, otherwise, the result is trivial.

\smallskip

Apply \thmref{thm:loj} to the function from $\U\times\U$ to $\R$ given by $g^2(x)+h^2(w)$. We conclude that there exists an open neighborhood $\W$ of $p$ in $\U$, and $c,\theta\in\R$ such that $c>0$, $0<\theta<1$, and, for all $(x, w)\in\W\times\W$,
$$
|g^2(x)+h^2(w)|^\theta\leq c\big|2g(x)(\nabla g(x), 0)+2h(w)(0,\nabla h(w)\big|.
$$
Restricting to the diagonal, we obtain that, for all $x\in\W$,
$$
|g^2(x)+h^2(x)|^\theta\leq 2c\sqrt{g^2(x)|\nabla g(x)|^2+h^2(x)|\nabla h(x)|^2}.
$$
Applying the Cauchy-Schwarz inequality, we have, for all $x\in\W$,
$$
|g^2(x)+h^2(x)|^\theta\leq 2c\sqrt{|(g^2(x), h^2(x))|\cdot |(|\nabla g(x)|^2, |\nabla h(x)|^2)|} = 2c\big(g^4(x)+h^4(x)\big)^{1/4}\big(|\nabla g(x)|^4+|\nabla h(x)|^4\big)^{1/4}.
$$

For $a, b, x>0$, $(a^x+b^x)^{1/x}$ is a decreasing function of $x$. Therefore, we conclude that, for all $x\in\W$,
$$
|(g(x), h(x))|^{2\theta}=|g^2(x)+h^2(x)|^\theta\leq 2c\big(g^2(x)+h^2(x)\big)^{1/2}\big(|\nabla g(x)|^2+|\nabla h(x)|^2\big)^{1/2},
$$
and so, for all $x\in\W$ such that $f(x)\neq 0$,
$$
|f(x)|^{2\theta-1}\leq 2c \sqrt{|\nabla g(x)|^2+|\nabla h(x)|^2}.
$$
This proves the result, except for the bounds on the exponent. As $-1<2\theta-1<1$, we have only to eliminate the possibility that $-1<2\theta-1\leq 0$. However, this is immediate, as we are assuming that $f(p)=0$ and  $|\nabla g(p)|^2+|\nabla h(p)|^2= 0$.
\end{proof}

\medskip

The following corollary follows at once.

\begin{cor}\label{cor:lojtuple2} Suppose that $G$ is real analytic and, for all $x\in \U$, $\nabla g_1(x), \dots, \nabla g_k(x)$ have the same magnitude.

Then, there exists an open neighborhood $\W$ of $p$ in $\U$, and $c,\theta\in\R$ such that $c>0$, $0<\theta<1$, and, for all $x\in\W$,
$$
|G(x)-G(p)|^\theta\leq c\big |\nabla g_1(x)\big|.
$$
\end{cor}

\bigskip

We now give the fundamental definition of this paper. Our intention is to isolate the properties of a complex analytic function that are used in proving the existence of the Milnor fibration inside a ball.

\begin{defn}\label{def:strongl} We say that $G$ satisfies the {\bf strong \L ojasiewicz inequality at $p$} or is an {\it \L-map at $p$}  if and only if 
there exists an open neighborhood $\W$ of $p$ in $\U$, and $c,\theta\in\R$ such that $c>0$, $0<\theta<1$, and, for all $x\in\W$,
$$
|G(x)-G(p)|^\theta\leq c\cdot\underset{|(a_1, \dots, a_k)|=1}{\operatorname{min}}\big |a_1\nabla g_1(x)+\dots+ a_k \nabla g_k(x)\big| = c\sigma_k(x).
$$

If $G$ satisfies the strong \L ojasiewicz inequality at $p$ and is real analytic in a neighborhood of $p$, then we say that {\bf $G$ is \L-analytic at $p$}.

We say that $G$ is {\bf  \L-analytic} if and only if $G$ is  \L-analytic at each point $p\in\U$.
\end{defn}

\medskip

\begin{rem}\label{rem:lanalrem} Note that if $\nabla g_1(x), \dots, \nabla g_k(x)$ are always pairwise orthogonal and have the same length, then the inequality in \defref{def:strongl} collapses to 
$$
|G(x)-G(p)|^\theta\leq c\big |\nabla g_1(x)\big|,
$$
which, as we saw in \corref{cor:lojtuple2}, is automatically satisfied if $G$ is real analytic.

In the case $k=2$, and $f=(g, h)$, one easily calculates the eigenvalues of the matrix
$$
\left[\begin{matrix}
|\nabla g|^2 & \nabla g\cdot\nabla h\\
\nabla g\cdot\nabla h & |\nabla h|^2
\end{matrix}\right ]
$$
and finds that $\sigma_2(x)$ is 
$$
\sqrt{\frac{|\nabla g(x)|^2+|\nabla h(x)|^2-\sqrt{(|\nabla g(x)|^2+|\nabla h(x)|^2)^2-4\big[|\nabla g(x)|^2|\nabla h(x)|^2-(\nabla g(x)\cdot \nabla h(x))^2\big]}}{2}},
$$
\end{rem}

\medskip

\begin{defn} When $G$ is such that $\nabla g_1(x), \dots, \nabla g_k(x)$ are always pairwise orthogonal and have the same length, we say that $G$ is a {\bf simple \L-map}.
\end{defn}

\smallskip

Thus, $G$ is a simple \L-map if and only if $M$ is a scalar multiple of the identity, i.e., for all $x\in\U$, $M(x)=|\nabla g_1(x)|^2 I_k$.

\smallskip

\begin{rem} Note that the pair of functions given by the real and imaginary parts of a holomorphic or anti-holomorphic function yield a simple \L-analytic map.
\end{rem}

\smallskip

\begin{defn}\label{def:slweight}  At each point $x\in\U$ where $d_xG\neq 0$ (i.e., where $|\nabla g_1(x)|^2+\dots+ |\nabla g_k(x)|^2\neq 0$), define
$$
\rho_G(x):=  \frac{k\left(\det M\right)^{1/k}}{\tr M}=\frac{k\left(\det M\right)^{1/k}}{|\nabla g_1(x)|^2+\dots+ |\nabla g_k(x)|^2}.
$$

\smallskip

If $d_xG$ is not identically zero, we define the {\bf \L-weight of $G$}, $\rho^{\operatorname{inf}}_G$, to be the infimum of  $\rho_G(x)$ over all $x\in\U$ such that $|\nabla g_1(x)|^2+\dots+ |\nabla g_k(x)|^2\neq 0$.

\smallskip

If $d_xG$ is identically zero, so that $G$ is locally constant, we set $\rho^{\operatorname{inf}}_G$ equal to $1$.

\smallskip

We say that {\bf $G$ has positive \L-weight at $x$} if and only if there exists an open neighborhood $\W$ of $x$ such that $G_{|_\W}$ has positive \L-weight.
\end{defn}

\smallskip

\begin{rem}\label{rem:lweight}  If $G$ is a submersion at $x$, then it is clear that $G$ has positive \L-weight at $x$. However, the converse is not true; for instance, a simple \L-map has positive \L-weight (see below). Thus, in a sense, having positive \L-weight is a generalization of being a submersion.

Note that, if $\rho^{\operatorname{inf}}_G>0$, we must have either $n\geq k$ or that $G$ is locally constant; furthermore, it must necessarily be true that $x\in\Sigma G$ implies that $\nabla g_1(x)=\cdots = \nabla g_k(x) =\0$. 

For a pair of functions $f=(g, h)$, one easily calculates that
$$
\rho_f(x)= \frac{2|\nabla g(x)||\nabla h(x)|\sin\eta(x)}{|\nabla g(x)|^2+ |\nabla h(x)|^2},
$$
where $\eta(x)$ equals the angle between $\nabla g(x)$ and $\nabla h(x)$.
\end{rem}

\medskip

One immediately concludes from \propref{prop:singval}:

\begin{prop} Suppose that $d_xG\neq 0$. Then, 
\begin{enumerate}
\item $0\leq \rho_G(x)\leq 1$;
\item $\rho_G(x)=0$ if and only if $\nabla g_1(x), \dots, \nabla g_k(x)$ are linearly dependent;
\item $\rho_G(x) =1$ if and only if $\nabla g_1(x), \dots, \nabla g_k(x)$ all have the same length and are pairwise orthogonal.
\end{enumerate}

Thus, $0\leq\rho^{\operatorname{inf}}_G\leq 1$, and $\rho^{\operatorname{inf}}_G=1$ if and only if $G$ is a simple \L-map.
\end{prop}

\medskip

In the case where $k=2$, there is a very precise relation between $\rho_f(x)$ and the strong \L-inequality.

\begin{prop} There exist positive constants $\alpha$ and $\beta$ such that, at all points $x\in\U$  such that  $d_xf\neq 0$,
$$
\alpha\sigma_2(x)\leq \rho_f(x)\sqrt{|\nabla g(x)|^2+|\nabla h(x)|^2}\leq \beta\sigma_2(x).
$$

Hence, $f$ satisfies the strong \L ojasiewicz inequality at $p$ if and only if 
there exists an open neighborhood $\W$ of $p$ in $\U$, and $c,\theta\in\R$ such that $c>0$, $0<\theta<1$, and, for all $x\in\W$ such that $d_xf\neq 0$,
$$
|f(x)-f(p)|^\theta\leq c\rho_f(x)\sqrt{|\nabla g(x)|^2+|\nabla h(x)|^2},
$$
and, for all $x\in\W$ such that $d_xf = 0$, $f(x)=f(p)$.
\end{prop}
\begin{proof} Assume that $d_xf\neq 0$. Then, the first set of inequalities is trivially true if $\sigma_2(x)=0$; so, assume $\sigma_2(x)\neq 0$. We need to show that we can find $\alpha, \beta>0$ such that
$$
\alpha \leq \frac{\rho_f(x)\sqrt{|\nabla g(x)|^2+|\nabla h(x)|^2}}{\sigma_2(x)}\leq\beta.
$$
This follows immediately from Item 2 of \propref{prop:singval} since
$$
\frac{\rho_f(x)\sqrt{|\nabla g(x)|^2+|\nabla h(x)|^2}}{\sigma_2(x)}= \frac{2\sqrt{\sigma_1^2(x)\sigma_2^2(x)}}{\sigma_2(x)\sqrt{|\nabla g(x)|^2+|\nabla h(x)|^2}}=\frac{2\sigma_1(x)}{\sqrt{|\nabla g(x)|^2+|\nabla h(x)|^2}}.
$$
\end{proof}

\medskip

\begin{prop}\label{prop:lmapequiv} The following are equivalent:
\begin{enumerate}
\item  $\rho^{\operatorname{inf}}_G>0$;
\item there exists $b>0$ such that, for all $x\in\U$, $b\cdot\sqrt{|\nabla g_1(x)|^2+\dots+ |\nabla g_k(x)|^2}\ \leq\ \sigma_k(x)$;
\item there exists $b>0$ such that, for all $x\in\U$, $b\cdot\sigma_1(x)\leq\sigma_k(x)$.
\end{enumerate}
\end{prop}
\begin{proof} Throughout this proof, when we take an infimum, we mean the infimum over all $x\in U$ such that  $d_xG\neq 0$; note that this implies that $\sigma_1(x)\neq 0$.  

Note that Items 2 and 3 simply say that
$$
0<\inf \frac{\sigma_k(x)}{\sqrt{|\nabla g_1(x)|^2+\dots+ |\nabla g_k(x)|^2}}
$$
and
$$
0<\inf \frac{\sigma_k(x)}{\sigma_1(x)}.
$$

Now,
$$\rho_G(x)= \frac{k\left(\det M\right)^{1/k}}{|\nabla g_1(x)|^2+\dots+ |\nabla g_k(x)|^2}=  \frac{k\left(\sigma_1^2(x)\cdots\sigma_k^2(x)\right)^{1/k}}{|\nabla g_1(x)|^2+\dots+ |\nabla g_k(x)|^2}=
$$
$$
k\left[\cdot \frac{\sigma_1(x)}{\sqrt{|\nabla g_1(x)|^2+\dots+ |\nabla g_k(x)|^2}}\cdots\frac{\sigma_k(x)}{\sqrt{|\nabla g_1(x)|^2+\dots+ |\nabla g_k(x)|^2}}\right]^{2/k},
$$
and the factors are non-negative, at most $1$, and are in decreasing order. This immediately yields the equivalence of Items 1 and 2.

The equivalence of Items 2 and 3 follows at once from
$$
\frac{\sigma_k(x)}{\sqrt{|\nabla g_1(x)|^2+\dots+ |\nabla g_k(x)|^2}}= \frac{\sigma_k(x)}{\sqrt{\sigma_1^2(x)+\cdots+\sigma_k^2(x)}}=
$$
$$ \frac{\sigma_k(x)/\sigma_1(x)}{\sqrt{1+ \left(\sigma_2(x)/\sigma_1(x)\right)^2+\cdots+\left(\sigma_k(x)/\sigma_1(x)\right)^2}}.
$$
\end{proof}

\medskip

Our primary interest in \L-weights is due to the following:

\begin{thm}\label{thm:posweight} Suppose that $G$ is real analytic and $\rho^{\operatorname{inf}}_G>0$. Then, $G$ is \L-analytic.
\end{thm}
\begin{proof} By \corref{cor:lojtuple}, if $x$ is near $p$ and $d_xG=0$, then $G(x)=G(p)$. Thus, it suffices to verify that the strong \L ojasiewicz inequality holds near any point at which $d_xG\neq 0$. The desired inequality follows at once from Item 2 of \propref{prop:lmapequiv}, combined with \corref{cor:lojtuple}.
\end{proof}

\medskip

Let $H:=(h_1, \dots, h_n)$ be a $C^1$ map from an open subset $\W$ of $\R^m$ to $\U$. For $x\in\W$, we let $A=A(H(x))$ be the matrix $[d_{H(x)}G]^t$, and set $B=B(x):=[d_xH]^t.$ Let $C=C(x)$ denote the matrix $[d_x(G\circ H)]^t = BA$.

\medskip

\begin{thm}\label{thm:lweightposcomp} Suppose that $q\in\W$, that $H$ which has positive \L-weight at $q$, and that $G$ has positive \L-weight at $H(q)$. Then, $G\circ H$ has positive \L-weight at $q$.
\end{thm}
\begin{proof} At a point $x\in\W$ where $C\neq 0$, \corref{cor:slw} tells us that
$$
\rho_{G\circ H}(x)\geq\frac{k\sigma_n^2(B(x))\sigma_k^2(A(H(x)))}{n\sigma_1^2(B(x))\sigma_1^2(A(H(x)))}.
$$
By Item 3 of \propref{prop:lmapequiv}, there exists an open neighborhood $\W^\prime$ of $q$ on which the infimum $b_H$ of $\sigma_n(B)/\sigma_1(B)$ is positive, and an open neighborhood $\U^\prime$ of $H(q)$ on which the infimum $b_G$ of $\sigma_n(A)/\sigma_1(A)$ is positive. Therefore, the infimum of $\rho_{G\circ H}(x)$ over all $x\in\W^\prime\cap H^{-1}(\U^\prime)$ such that $C(x)\neq 0$ is at least $kb_H^2b_G^2/n$.
\end{proof}

\medskip

\thmref{thm:lweightposcomp} gives us an easy way of producing new non-simple \L-analytic maps: take a map into $\R^2$ consisting of the real and imaginary parts of a holomorphic or anti-holomorphic function, and then compose with a real analytic change of coordinates on an open set in $\R^2$. Actually, we can give a much more precise result when $H$ is a simple \L-map.

\begin{prop}\label{prop:easycomp} Suppose that $H$ is a simple \L-map. Then, for all $x\in\W$ at which $d_x(G\circ H)\neq 0$, $\rho_{G\circ H}(x) = \rho_G(H(x))$.
\end{prop}
\begin{proof} Suppose that $d_x(G\circ H)\neq 0$. As $H$ is a simple \L-map, $B^tB=\lambda\, I_n$, where $\lambda= |\nabla h_1(x)|^2$. Thus,
$$
\rho_{G\circ H}(x)= \frac{k\left[\det(A^tB^tBA)\right]^{1/k}}{\tr(A^tB^tBA)}= \frac{k\left[\lambda^k\det(A^tA)\right]^{1/k}}{\lambda\,\tr(A^tA)} = \frac{k\left[\det(A^tA)\right]^{1/k}}{\tr(A^tA)} =\rho_G(H(x)).
$$
\end{proof}

\medskip

We give two quick examples of how to produce interesting \L-analytic maps.

\begin{exm}\label{exm:mix}  As we mentioned above, if $g$ and $h$ are the real and imaginary parts of a holomorphic or anti-holomorphic function, then $f=(g, h)$ is a simple \L-analytic function.

One can also mix holomorphic and anti-holomorphic functions. Let $z=x+iy$ and $w=u+iv$, and consider the real and imaginary parts of $\overline{z}w^2$, i.e., let $g= x(u^2-v^2)+2yuv$ and $h= 2xuv-y(u^2-v^2)$. It is trivial to verify that $(g, h)$ is simple \L-analytic.

More generally, if $(g_1, h_1)$ and $(g_2, h_2)$ are simple \L-analytic, then the real and imaginary parts of 
$$(g_1(\mathbf z)+ih_1(\mathbf z))(g_2(\mathbf w)+ih_2(\mathbf w))
$$
will yield a new simple \L-analytic function.
\end{exm}

\begin{exm} Suppose that $f=(g, h)$ is a simple \L-analytic map. 

Consider the linear map $L$ of $\R^2$ given by $(u, v)\mapsto (au+bv, cu+dv)$. Then, at every point, the matrix of the derivative of $L$ is $\left[\begin{matrix}a & b\\c & d\end{matrix}\right]$, and one calculates that $\rho_L$ is constantly $2|ad-bc|/(a^2+b^2+c^2+d^2)$. Thus, if $ad-bc\neq 0$, then $\rho_L>0$  and, unless $ab+cd=0$ and $a^2+c^2 = b^2+d^2$, we have that $\rho_L<1$.

Therefore, if $L$ is an isomorphism which is not an orthogonal transformation composed with scalar multiplication, then \propref{prop:easycomp} tells us that $(ag+bh, cg+dh)$ has positive \L-weight less than one and, hence, is an \L-analytic map which is not simple and, therefore, cannot arise from a holomorphic or anti-holomorphic complex function.
\end{exm}

\bigskip

\noindent{\bf Relationship to Jacquemard's Conditions}

\smallskip

For the remainder of this section, we will restrict ourselves to considering a real analytic map $f=(g,h)$ into $\R^2$. In \cite{jacquemard}, Jacquemard investigates such $f$ satisfying three conditions. We will refer to these conditions as $J0$, $J1$ and $J2$; actually, our condition $J2$ will be the weaker condition given by Ruas, Seade, and Verjovsky in \cite{ruasseadeverj}.

\begin{defn}\label{def:jconditions} The {\bf Jacquemard conditions} are:

\noindent $(J0)$ the origin is an isolated critical point of $f$;

\noindent $(J1)$ there exists an open neighborhood $\W$ of $0$ in $\U$ and a real number $\tau>o$ such that, if $x\in\W$ is such that $\nabla g(x)\neq\mathbf 0$ and $\nabla h(x)\neq\mathbf 0$, then
$$
\left|\frac{\nabla g(x)}{|\nabla g(x)|}\cdot\frac{\nabla h(x)}{|\nabla h(x)|}\right|\leq 1-\tau;
$$

\noindent $(J2)$ the real integral closures of the Jacobian ideals of $g$ and $h$, inside the ring of real analytic germs at the origin, are equal.
\end{defn}

\medskip

The reader should note that we may discuss the Jacquemard conditions holding independently; in particular, we will assume in some settings that $J1$ and $J2$ hold, without assuming $J0$.

\medskip

The importance of the Jacquemard conditions stems from the following theorem, which is essentially proved in \cite{jacquemard}, and improved using J2 in \cite{ruasseadeverj}.

\begin{thm}  If Jacquemard's conditions hold, then $f$ satisfies the {\it strong Milnor condition}, i.e., there exists $\epsilon_0>0$ such that, for all $\epsilon$ such that $0<\epsilon\leq\epsilon_0$, the map $f/|f|$ from $S_\epsilon$ (the sphere of radius $\epsilon$, centered at the origin) to $\R^2$ is a smooth, locally trivial fibration. 
\end{thm}

\medskip

We wish to see that, if we assume $J1$, greatly weaken $J2$, and omit $J0$, then $f$ has positive \L-weight at $0$ and, hence, is \L-analytic at $0$. It will then follow from \corref{cor:main} that the Milnor fibration inside a ball, centered at $0$, exists. First, we need a lemma.

\begin{lem}\label{lem:realintclos} Condition J2 implies that there exists an open neighborhood $\W$ of the origin in $\U$ in which $\nabla g$ and $\nabla h$ are comparable in magnitude, i.e., such that there exist $A, B >0$ such that, for all $x\in\W$, 
$$A|\nabla g(x)|\leq |\nabla h(x)|\leq B|\nabla g(x)|.$$
\end{lem}
\begin{proof} By Proposition 4.2 of \cite{gaffinvent}, J2 is equivalent to: there exists a neighborhood $\W$ of the origin in $\U$ and $C_1, C_2>0$ such that, at all points in $\W$, for all $i$ such that $1\leq i\leq n$, $|\partial g/\partial x_i|\leq C_1\max_j|\partial h/\partial x_j|$ and $|\partial h/\partial x_i|\leq C_2\max_j|\partial g/\partial x_j|$. One quickly concludes that, at all points of $\W$, $|\nabla g|\leq\sqrt{n}C_1|\nabla h|$ and $|\nabla h|\leq\sqrt{n}C_2|\nabla g|$. The lemma follows.
\end{proof}

\smallskip

\begin{prop}\label{prop:jacequiv} Suppose there exists an open neighborhood $\W$ of the origin in $\U$ in which $\nabla g$ and $\nabla h$ are comparable in magnitude. Then, Condition J1 is satisfied if and only if $f$ has positive \L-weight at $0$.
\end{prop}
\begin{proof} First, note that, when $\nabla g$ and $\nabla h$ are comparable in magnitude, then one of them equals zero at a point $x$ if and only if both of them equal  zero at $x$. Suppose that $A, B >0$ are such that, for all $x\in\W$, 
$$A|\nabla g(x)|\leq |\nabla h(x)|\leq B|\nabla g(x)|.$$

Let $\eta(x)$ equal the angle between $\nabla g(x)$ and $\nabla h(x)$. Then, 
$$\left|\frac{\nabla g(x)}{|\nabla g(x)|}\cdot\frac{\nabla h(x)}{|\nabla h(x)|}\right| = |\cos\eta(x)|,$$
and it is trivial to conclude that Condition J1 holds if and only if there exists an open neighborhood of the origin $\W^\prime\subseteq\W$ such that $0<\inf_{x\in\W^\prime}\big(\sin\eta(x)\big)$.

Recall from \remref{rem:lweight} that, at a point $x$ where $|\nabla g(x)|^2+ |\nabla h(x)|^2\neq 0$,
$$
\rho_f(x)= \frac{2|\nabla g(x)||\nabla h(x)|\sin\eta(x)}{|\nabla g(x)|^2+ |\nabla h(x)|^2}.
$$
Now, $A|\nabla g(x)|\leq |\nabla h(x)|\leq B|\nabla g(x)|$ implies that
$$
\frac{2A}{1+B^2}\leq \frac{2|\nabla g(x)||\nabla h(x)|}{|\nabla g(x)|^2+ |\nabla h(x)|^2}\leq \frac{2B}{1+A^2}.
$$
The conclusion follows.
\end{proof}

\begin{rem} In light of \lemref{lem:realintclos} and \propref{prop:jacequiv}, we see that the Jacquemard conditions are far stronger than what one needs to conclude that $f$ has positive \L-weight at the origin.
\end{rem}

\section{Milnor's conditions (a) and (b)}\label{sec:milnorcondab}

In this section, we will discuss general conditions which allow us to conclude that Milnor fibrations exist.

\smallskip

We continue with our notation from the previous section, except that we now assume that $G$ is only $C^\infty$. For notational convenience, we restrict our attention to the case where the point $p$ is the origin and where $G(\0)=\0$. We also assume that $G$ is not locally constant near $\0$. Let $X:=G^{-1}(\0)=V(G)$.

Let $\mathfrak A=\Sigma G$, so that $\mathfrak A$ is the closed set of points in $\U$ at which the gradients $\nabla g_1, \nabla g_2, \dots \nabla g_k$ are linearly dependent.  Let $r$ denote the function given by the square of the distance from the origin, and let $\mathfrak B$ denote the closed set of points in $\U$ at which the gradients $\nabla r, \nabla g_1, \nabla g_2, \dots \nabla g_k$ are linearly dependent. Of course, $\mathfrak A\subseteq \mathfrak B$. 

\smallskip

We wish to give names to the submersive conditions necessary to apply Ehresmann's Theorem \cite{ehresmann} (in the case of manifolds with boundary).

\begin{defn}\label{def:milnorab} We say that the map $G$ satisfies {\bf Milnor's condition (a) at $\0$} (or that {\bf $\0$ is an isolated critical value of $G$ near $\0$}) if and only if $\0\not\in\overline{\mathfrak A-X}$, i.e., if $\Sigma G\subseteq V(G)$ near $\0$.

We say that the map $G$ satisfies {\bf Milnor's condition (b) at $\0$}  if and only if $\0$ is an isolated point of (or, is not in) $X\cap\overline{\mathfrak B-X}$.

If $G$ satisfies Milnor's condition (a) (respectively,  (b)), then we say that $\epsilon>0$ is a {\bf Milnor (a) (respectively, (b)) radius for $G$ at $\0$} provided that $B_\epsilon\cap (\overline{\mathfrak A-X}) =\emptyset$ (respectively, $B_\epsilon\cap X\cap(\overline{\mathfrak B-X})\subseteq\{\0\}$). We say simply that {\bf $\epsilon>0$ is a Milnor radius for $G$ at $\0$} if and only if $\epsilon$ is both a Milnor (a) and Milnor (b) radius for $G$ at $\0$.
\end{defn}

\begin{rem} Using our notation from the introduction, if $f_{{}_\C}$ is complex analytic, then $f=(g,h)$ satisfies Milnor's conditions (a) and (b); in this case, Milnor's condition (a) is well-known and follows easily from a curve selection argument or from \corref{cor:complexloj}, and Milnor's condition (b) follows from the existence of good (or $a_f$) stratifications of $V(f)$ (see \cite{hammlezariski} and \cite{relmono}, and below).
\end{rem}

\smallskip

\begin{lem}\label{lem:milfib}
Suppose that the map $G$ satisfies Milnor's condition (b) at $\0$ and that $\epsilon_0>0$ is a Milnor (b) radius for $G$ at $\0$.  Then, for all $\epsilon$  such that $0<\epsilon<\epsilon_0$, there exists $\delta_\epsilon>0$ such that the map
$$H: (\ob_{\epsilon_0}-B_{\epsilon})\cap G^{-1}\big(\ob_{\delta_\epsilon}-\{\0\}\big)\rightarrow(\ob_{\delta_\epsilon}-\{\0\})\times (\epsilon^2, \epsilon_0^2)$$
given by $H(x)=(G(x), |x|^2)$ is a proper submersion.

In particular, for all $\epsilon^\prime$ such that $0<\epsilon^\prime<\epsilon_0$, there exists $\delta_{\epsilon^\prime}>0$ such that 
$$G: \partial B_{\epsilon^\prime}\cap G^{-1}\big(\ob_{\delta_{\epsilon^\prime}}-\{\0\}\big)\rightarrow \ob_{\delta_{\epsilon^\prime}}-\{\0\}
$$
is a proper submersion.
\end{lem}
\begin{proof} That $H$ is proper is easy. Let $\pi: (\ob_{\delta_\epsilon}-\{\0\})\times (\epsilon^2, \epsilon_0^2)\rightarrow (\epsilon^2, \epsilon_0^2)$ denote the projection. Suppose that $C\subseteq (\ob_{\delta_\epsilon}-\{\0\})\times (\epsilon^2, \epsilon_0^2)$ is compact. Then, $\pi(C)$ is compact, and $H^{-1}(C)$ is a closed subset of the compact set $\{x\in \ob_{\epsilon_0}\ |\ |x|^2\in\pi(C)\}$. Thus, $H^{-1}(C)$ is compact.

Now, a critical point of $H$ is precisely a point in $\mathfrak B\cap  (\ob_{\epsilon_0}-B_{\epsilon})\cap G^{-1}\big((\ob_{\delta_\epsilon}-\{\0\})\times (\epsilon^2, \epsilon_0^2)\big)$. Suppose that had such a point, regardless of how small we choose $\delta_\epsilon>0$. Then, we would have a sequence $x_i\in (\mathfrak B-X)\cap (\ob_{\epsilon_0}-B_{\epsilon})$ such that $G(x_i)\rightarrow 0$. As the $x_i$ are in the compact set $B_{\epsilon_0}-\ob_{\epsilon}$, by taking a subsequence if necessary, we may assume that $x_i\rightarrow x\in B_{\epsilon_0}-\ob_{\epsilon}$. As $G(x_i)\rightarrow 0$, $x\in X$. Therefore, $x\in (B_{\epsilon_0}-\ob_{\epsilon})\cap X\cap\overline{(\mathfrak B-X)}$; a contradiction, since $\epsilon_0$ is a Milnor (b) radius. Hence, $H$ is a submersion.

The last statement follows at once from this, since one need only pick an $\epsilon$ such that $0<\epsilon<\epsilon^\prime$ and apply that $H$ is a submersion.
\end{proof}

\smallskip

\begin{thm}\label{thm:milfibab} Suppose that $G$ satisfies Milnor's conditions (a) and (b) at $\0$, and let $\epsilon_0$ be a Milnor radius for $G$ at $\0$.

Then, for all $\epsilon$  such that $0<\epsilon<\epsilon_0$, there exists $\delta_\epsilon>0$ such that the map
$$
G:B_{\epsilon}\cap G^{-1}\big(\ob_{\delta_{\epsilon}}-\{\0\}\big)\rightarrow \ob_{\delta_\epsilon}-\{\0\}
$$
is a proper, stratified submersion and, hence, a locally trivial fibration, in which the local trivializations preserve the strata.

In addition, for all such $(\epsilon, \delta_\epsilon)$ pairs, for all $\delta$ such that $0<\delta<\delta_\epsilon$, the map
$$
G:B_{\epsilon}\cap G^{-1}\big(\partial B_\delta\big)\rightarrow \partial B_{\delta}
$$
is a proper, stratified submersion and, hence, a locally trivial fibration, whose diffeomorphism-type is independent  of the choice of such $\epsilon$ and $\delta$.

It follows that, for all such $(\epsilon, \delta_\epsilon)$ pairs, for all $\delta$ such that $0<\delta<\delta_\epsilon$, the map
$$
G:\ob_{\epsilon}\cap G^{-1}\big(\partial B_\delta\big)\rightarrow \partial B_{\delta}
$$
is a locally trivial fibration, whose diffeomorphism-type is independent  of the choice of such $\epsilon$ and $\delta$.
\end{thm}
\begin{proof} The map $
G:B_{\epsilon}\cap G^{-1}\big(\ob_{\delta_{\epsilon}}-\{\0\}\big)\rightarrow \ob_{\delta_\epsilon}-\{\0\}$ is clearly proper.

Milnor's condition (a) tells us that
$$
G:\ob_{\epsilon}\cap G^{-1}\big(\ob_{\delta_{\epsilon}}-\{\0\}\big)\rightarrow \ob_{\delta_\epsilon}-\{\0\},
$$
is a submersion regardless of the choice of $\delta_{\epsilon}>0$.

The last line of \lemref{lem:milfib} tells us that we may pick $\delta_\epsilon>0$ such that the map
$$
G:\partial B_{\epsilon}\cap G^{-1}\big(\ob_{\delta_{\epsilon}}-\{\0\}\big)\rightarrow \ob_{\delta_\epsilon}-\{\0\}
$$
is a submersion.

Therefore,
$$
G:B_{\epsilon}\cap G^{-1}\big(\ob_{\delta_{\epsilon}}-\{\0\}\big)\rightarrow \ob_{\delta_\epsilon}-\{\0\}
$$
is a proper, stratified submersion, and so, by Ehresmann's Theorem (with boundary) \cite{ehresmann} or Thom's first isotopy lemma \cite{mather}, this map is a locally trivial fibration, in which the local trivializations preserve the strata.

\smallskip

It follows at once that for a fixed such $\epsilon$, for all $\delta$ such that $0<\delta<\delta_\epsilon$, the map
$$
G:B_{\epsilon}\cap G^{-1}\big(\partial B_\delta\big)\rightarrow \partial B_{\delta}
$$
is a proper, stratified submersion and, hence, a locally trivial fibration, whose diffeomorphism-type is independent  of the choice of such $\delta$.

\smallskip

It remains for us to show that, if we pick $0<\epsilon<\epsilon^\prime<\epsilon_0$, then there exists $\delta>0$ such that $\delta<\min\{\delta_\epsilon, \delta_{\epsilon^\prime}\}$, and such that the fibrations $
G:B_{\epsilon}\cap G^{-1}\big(\partial B_\delta\big)\rightarrow \partial B_{\delta}$ and $
G:B_{\epsilon^\prime}\cap G^{-1}\big(\partial B_\delta\big)\rightarrow \partial B_{\delta}
$ have the same diffeomorphism-type.

Let $\hat\epsilon$ be such that $0<\hat\epsilon<\epsilon$, let $\delta_{\hat\epsilon}>0$ be as in the first part of \lemref{lem:milfib}, and let $\delta>0$ be less than $\min\{\delta_\epsilon, \delta_{\epsilon^\prime}, \delta_{\hat\epsilon}\}$. Then, as the interval $(\hat\epsilon^2, \epsilon_0^2)$ is contractible, \lemref{lem:milfib} tells us immediately that $
G:B_{\epsilon}\cap G^{-1}\big(\partial B_\delta\big)\rightarrow \partial B_{\delta}
$ and $
G:B_{\epsilon^\prime}\cap G^{-1}\big(\partial B_\delta\big)\rightarrow \partial B_{\delta}
$ have the same diffeomorphism-type.
\end{proof}

\smallskip

\begin{rem} One might hope suspect that, if $n$ and $k$ are even, then Milnor's condition's (a) and (b) are satisfied by maps $G$ which come from complex analytic maps. This is {\bf not} the case.

Even in the nice case of a complex analytic isolated complete intersection singularity, the set of critical values would not locally consist solely of the origin, but would instead be a hypersurface in an open neighborhood of the origin in $\C^{k/2}$; see \cite{looibook}. Thus, Milnor's condition (a) does not hold.

This means that the types of Milnor fibrations that we obtain when we assume Milnor's conditions (a) and (b) are extremely special.
\end{rem}

\smallskip

\begin{defn} If $G$ satisfies Milnor's conditions (a) and (b) at $\0$, we refer to the two fibrations over spheres (or their diffeomorphism classes) in \thmref{thm:milfibab}, $
G:B_{\epsilon}\cap G^{-1}\big(\partial B_\delta\big)\rightarrow \partial B_{\delta}
$ and $
G:\ob_{\epsilon}\cap G^{-1}\big(\partial B_\delta\big)\rightarrow \partial B_{\delta}$,
as the {\bf compact Milnor fibration of $G$ at $\0$ inside a ball} and the {\bf Milnor fibration of $G$ at $\0$ (inside a ball)}, respectively.
\end{defn}

\smallskip

\begin{cor}\label{cor:open} Suppose that $G$ satisfies Milnor's conditions (a) and (b) at $\0$, and that $k>1$. Then, $G$ maps an open neighborhood of the origin onto an open neighborhood of the origin.\end{cor}
\begin{proof}  Recall that $G(\0)=\0$ and that $G$ is not locally constant by assumption. As $k>1$, $\ob_{\delta_\epsilon}-\{\0\}$ is connected, and \thmref{thm:milfibab} implies that $G: \ob_\epsilon\cap G^{-1}(\ob_{\delta_\epsilon})\rightarrow \ob_{\delta_\epsilon}$ is surjective.
\end{proof}

\begin{defn} Suppose that $G$ satisfies Milnor's conditions (a) and (b) at $\0$, and that $\epsilon$ is a Milnor radius for $G$ at $\0$. If $\delta>0$, then $(\epsilon, \delta)$ is a {\bf Milnor pair for $G$ at $\0$} if and only if there exists a $\hat\delta>\delta$ such that $G: B_\epsilon\cap G^{-1}(\ob_{\hat\delta}^*)\rightarrow \ob_{\hat\delta}^*$ is a stratified submersion (which is, of course, smooth and proper).
\end{defn}

\begin{rem} Whenever we write that $(\epsilon, \delta)$ is a Milnor pair for $G$ at $\0$, we are assuming that $G$ satisfies Milnor's conditions (a) and (b) at $\0$.
\end{rem}

\section{The Main Theorem}\label{sec:main}

Now that we have finished our general discussion of \L-analytic functions and Milnor's conditions, we wish to investigate how they are related to each other.

\medskip

The following proposition is trivial to conclude.

\begin{prop}\label{prop:milnora} If $G$ is a $C^1$ \L-map at $p$, then, near $p$, $\Sigma G\subseteq G^{-1}(G(p))$, i.e., $G$ satisfies Milnor's condition (a) at $p$. 
\end{prop}

\medskip

Our goal is to prove that if $G$ is \L-analytic at $p$, then $G$ also satisfies Milnor's condition (b) at $p$, for then \thmref{thm:milfibab} will tell us that the Milnor fibrations inside a ball exist.

It will turn out that our method of proof allows us to conclude this for $k=1$, $2$, $4$, and $8$, by using the normed division algebra structures in those dimensions. In fact, what we need is the following:

\begin{defn}\label{def:napmap} Let $\W$ be an open neighborhood of the origin in $\R^k$. Let $P$ be an infinite subset of $\N$ (the {\bf admissible powers}). A function $M:P\times\W\rightarrow\R^k$, where $M(p, y)$ is written as $M_p(y)$, is a {\bf normed analytic power map} (a {\bf nap-map}) if and only if, for all $p\in P$, \begin{enumerate}
\item $M_p$ is real analytic;
\item for all $y\in \W$, $|M_p(y)|=|y|^p$;
\item the image of $M_p$ contains an open neighborhood of the origin;
and
\item there exists a constant $K_p>0$ such that, if $M_p=(m_1, \dots, m_k)$, then, for all $y\in \W$,
$$
\underset{|(t_1, \dots, t_k)|=1}{\operatorname{max}}\big |t_1\nabla m_1(y)+\dots +t_k\nabla m_k(y)\big|\leq K_p|y|^{p-1}.
$$
\end{enumerate}
\end{defn}

\medskip

\begin{rem}\label{rem:minmax} It is easy to show that Items 1 and 2 in the above definition (or even replacing real analytic by $C^1$) imply that
$$
\underset{|(t_1, \dots, t_k)|=1}{\operatorname{min}}\big |t_1\nabla m_1(y)+\dots +t_k\nabla m_k(y)\big|\leq\ p|y|^{p-1}\ \leq\underset{|(t_1, \dots, t_k)|=1}{\operatorname{max}}\big |t_1\nabla m_1(y)+\dots +t_k\nabla m_k(y)\big|.
$$
We will not use these inequalities.
\end{rem}
\medskip

Originally, we used multiplication in the reals, complex numbers, quaternions, and octonions to produce nap-maps from $\R^k$ to $\R^k$ when $k=1$, $2$, $4$, or $8$. T. Gaffney produced the easier nap-maps below {\bf for all $k$}.

\medskip

\begin{prop}\label{prop:napmap} For all $k$, there exists a nap-map from $\R^k$ to $\R^k$. In particular, if $P$ is the set of odd natural numbers, then, for all $k$, the function $M:P\times \R^k\rightarrow\R^k$ given by $M(p, y)=M_p(y) = |y|^{p-1}y$ is a nap-map.
\end{prop}
\begin{proof} Let $P$ and $M_p$ be as in the statement of the proposition. Since $p\in P$ is odd, Items 1, 2, and 3 in the definition of a nap-map are trivially satisfied. We need to prove the inequality in Item 4 of \defref{def:napmap}.

Let $M_p=(m_1, \dots, m_k)$. One calculates that
$$
\frac{\partial m_i}{\partial y_j} = |y|^{p-1}\delta_{i,j} +(p-1)|y|^{p-3}y_iy_j,
$$
where $\delta_{i,j}$ is the Kronecker delta function. Thus, considering $y$ as a column-vector and denoting its transpose by $y^t$, the transpose of the derivative matrix of $M_p$ at $y$ is
$$
L:= |y|^{p-1}I +(p-1)|y|^{p-3}y^ty.
$$

Now, suppose that $T$ is a column vector of unit length in $\R^k$. Our goal is to show that $|LT|\leq K_p|y|^{p-1}$ for some constant $K_p$. We find
$$
|LT|=\big|\,|y|^{p-1}T+(p-1)|y|^{p-3}y^tyT\big|\leq \big|\,|y|^{p-1}T\big|+\big|(p-1)|y|^{p-3}y^tyT\big|\leq |y|^{p-1}+ (p-1)|y|^{p-1},
$$
where we use repeatedly in the last inequality that for matrices $A$ and $B$, $|AB|\leq |A|\,|B|$.
Thus,
$\underset{|(t_1, \dots, t_k)|=1}{\operatorname{max}}\big |t_1\nabla m_1(y)+\dots +t_k\nabla m_k(y)\big| = |L|\leq p|y|^{p-1}$.
\end{proof}

\bigskip

We can now prove our main lemma.

\begin{lem}\label{lem:main} Suppose that $G$ is \L-analytic at $\0$, and $G(\0)=\0$. Then, exists an open neighborhood $\W$ of $\0$ in $\U$ and a Whitney stratification $\strat$ of $\W\cap X$ such that, for all $S\in\strat$, the pair $(\W-X, S)$ satisfies Thom's $a_G$ condition, i.e., if $p_j\rightarrow p\in S\in\strat$, where $p_j\in\W-X$, and $T_{p_j}G^{-1}G(p_j)$ converges to some linear subspace $\cal T$ in the Grassmanian of $(n-k)$-dimensional linear subspaces of $\R^n$, then $T_pS\subseteq\cal T$.
\end{lem}
\begin{proof} This proof follows that of Theorem 1.2.1 of  \cite{hammlezariski}. 

\smallskip

Let $\W$, $c$, and $\theta$ be as in \defref{def:strongl}: $\W$ is an open neighborhood of $\0$ in $\U$, and $c,\theta\in\R$ are such that $c>0$, $0<\theta<1$, and, for all $x\in\W$,
$$
|G(x))|^\theta\leq c\cdot\underset{|(a_1, \dots, a_k)|=1}{\operatorname{min}}\big |a_1\nabla g_1(x)+\dots+ a_k \nabla g_k(x)\big| = c\sigma_k(x).
$$
Let $M$  be a nap-map on a neighborhood $\V$ of the origin in $\R^k$ (which exist by \propref{prop:napmap}), and let $P\subseteq\N$ denote the set of admissible powers. Let $\pi\in P$ be such that $\pi>1/(1-\theta)$, so that the image of $M_\pi$ contains a neighborhood of the origin. Let $M_\pi=(m_1, \dots, m_k)$.

Consider the real analytic map $\widehat G:\W\times\V\rightarrow\R^k$ given by 
$$
\widehat G(x,y):= (g_1(x)+m_1(y), \dots, g_k(x)+m_k(y)).
$$

Let $\widehat\strat$ be a Whitney stratification of $V(\widehat G)$ such that $V(G)\times\{\0\}$ is a union of strata. Let $\strat:=\{S\ |\ S\times\{\0\}\in\widehat\strat\}$. We claim that $\W$ and $\strat$ satisfy the conclusion of the lemma.

It will be convenient to deal with the conormal formulation of the $a_G$ condition. Suppose that we have a sequence of points $p_j\rightarrow p\in S\in\strat$, where $p_j\in\W-X$, and a sequence ${}^j\hskip -.01in \mathbf a:=({}^j\hskip -.01in a_1, \dots {}^j\hskip -.01in a_k)\in\R^k$ such that ${}^j\hskip -.01in a_1d_{p_j}g_1+\dots +{}^j\hskip -.01in a_kd_{p_j}g_k$ converges to a cotangent vector $\eta$ (in the fiber of $T^*\W$ over $p$). We wish to show that $\eta(T_pS)\equiv 0$. 

\smallskip

If $\eta\equiv 0$, there is nothing to show; so we assume that $\eta\not\equiv 0$. 

\smallskip

We may assume that, for all $j$, ${}^j\hskip -.01in \mathbf a\neq \0$. Also, as $p_j\rightarrow p\in V(G)$, we may assume that, for all $j$, $-G(p_j)$ is in the image of $M_\pi$. Let ${}^j\hskip -.01in \mathbf u:=({}^j\hskip -.01in u_1, \dots, {}^j\hskip -.01in u_k)\in M_\pi^{-1}(-G(p_j))$, so that $q_j:=(p_j, {}^j\hskip -.01in \mathbf u)\in V(\widehat G)$. Since $G(p)=0$ and $\widehat G(q_j)=0$, $M_\pi({}^j\hskip -.01in \mathbf u)\rightarrow \0$; by Item 2 in the definition of a nap-map, it follows that ${}^j\hskip -.01in \mathbf u\rightarrow\0$, and so $q_j\rightarrow (p, \0)\in S\times\{\0\}\in\widehat\strat$.

By taking a subsequence, if necessary, and using that $\eta\not\equiv 0$, we may assume that the projective class 
$$
\left[{}^j\hskip -.01in a_1d_{q_j}(g_1+m_1)+\dots+{}^j\hskip -.01in a_kd_{q_j}(g_k+m_k)\right]
$$
converges to some $[\omega]$, where $\omega= b_1dx_1+\dots b_ndx_n+ s_1dy_1+\dots+s_kdy_k\not\equiv 0$. By Whitney's condition (a), we have that $\omega(T_pS\times \{\0\})\equiv 0$. 

\smallskip

Let us reformulate part of our discussion above using vectors, instead of covectors. In terms of vectors, we are assuming that ${}^j\hskip -.01in a_1\nabla g_1(p_j)+\dots +{}^j\hskip -.01in a_k\nabla g_k(p_j)\rightarrow\mathbf w\neq\0$, and that
$$
\left[{}^j\hskip -.01in a_1\Big(\nabla g_1(p_j), \nabla m_1({}^j\hskip -.01in \mathbf u)\Big)+\dots+ {}^j\hskip -.01in a_k\Big(\nabla g_k(p_j), \nabla m_k({}^j\hskip -.01in \mathbf u)\Big)\right]\rightarrow \left[(b_1, \dots, b_n, s_1, \dots, s_k)\right]\neq \0.
$$
If we could show that the projective classes $[\mathbf w\times\{\0\}]$ and $\left[(b_1, \dots, b_n, s_1, \dots, s_k)\right]$ are equal, then we would be finished.

To show that $[\mathbf w\times\{\0\}]$ and $\left[(b_1, \dots, b_n, s_1, \dots, s_k)\right]$ are equal, it clearly suffices to show that
$$
\frac
{\left|{}^j\hskip -.01in a_1\nabla m_1({}^j\hskip -.01in \mathbf u)+\dots+ {}^j\hskip -.01in a_k\nabla m_k({}^j\hskip -.01in \mathbf u)\right|}
{\left|{}^j\hskip -.01in a_1\nabla g_1(p_j)+\dots+ {}^j\hskip -.01in a_k\nabla g_k(p_j)\right|}\rightarrow 0.
$$
Dividing the numerator and denominator by $|{}^j\hskip -.01in \mathbf a|$, we see that we may assume that $|{}^j\hskip -.01in \mathbf a|=1$.

Using that $G$ is an \L-map at $\0$ and Item 4 in the definition of a nap-map, we obtain that
$$
\frac
{\left|{}^j\hskip -.01in a_1\nabla m_1({}^j\hskip -.01in \mathbf u)+\dots+ {}^j\hskip -.01in a_k\nabla m_k({}^j\hskip -.01in \mathbf u)\right|}
{\left|{}^j\hskip -.01in a_1\nabla g_1(p_j)+\dots+ {}^j\hskip -.01in a_k\nabla g_k(p_j)\right|}\ \leq\ \frac{K_p|{}^j\hskip -.01in \mathbf u|^{\pi-1}}{|G(p_j)|^\theta}.
$$

Now, $|{}^j\hskip -.01in \mathbf u|^{\pi}=|M_\pi({}^j\hskip -.01in \mathbf u)| = |G(p_j)|$. Thus, $|{}^j\hskip -.01in \mathbf u|^{\pi-1} = |G(p_j)|^{(\pi-1)/\pi}$, and we would like to show that 
$$
\frac{ |G(p_j)|^{(\pi-1)/\pi}}{|G(p_j)|^\theta}\rightarrow 0.
$$
However, this follows at once, since $G(p_j)\rightarrow\0$ and, by our choice of $\pi$, $(\pi-1)/\pi >\theta$.
\end{proof}

\medskip

\begin{rem}\label{rem:subanal} In \lemref{lem:main}, we assumed that $G$ was real analytic and used that a normed, analytic power map $M_p$ exists on $\R^k$. We assumed analyticity for both maps so that $V(\widehat G)$ would have a Whitney stratification. However, it is enough to assume that $G$ and $M_p$ are {\it subanalytic}. See \cite{hardtsubanal}.
\end{rem}

\medskip

The main theorem follows easily. Note that, when $k=1$, \thmref{thm:loj} implies that, if $G$ is real analytic, then $G$ is \L-analytic.

\begin{thm}\label{thm:main} Suppose that $G$ is \L-analytic at $p$. Then, Milnor's condition's (a) and (b) hold at $p$.
\end{thm}

\begin{proof} We assume, without loss of generality, that $p=\0$ and $G(\0)=\0$.

Milnor's condition (a) is immediate, as we stated in \propref{prop:milnora}. Milnor's condition (b) follows from \lemref{lem:main} by using precisely the argument of L\^e in \cite{relmono}. 

Let $\epsilon_0>0$ be a Milnor (a) radius for $G$ at $\0$ such that, for all $\epsilon^\prime$ such that $0<\epsilon^\prime\leq\epsilon_0$, the sphere $S_{\epsilon^\prime}$ transversely intersects all of the strata of the Whitney $a_G$ stratification whose existence is guaranteed by \lemref{lem:main}. We claim that $\epsilon_0$ is a Milnor (b) radius. 

Suppose not. Then, there would be a sequence of points $p_i\in B_{\epsilon_0}-V(G)$ such that $G(p_i)\rightarrow 0$ and $T_{p_i}\big(G^{-1}G(p_i)\big)\subseteq T_{p_i}S_{\epsilon_i}$, where $\epsilon_i$ denotes the distance of $p_i$ from the origin (we use Milnor's condition (a) here to know that $G^{-1}G(p_i)$ is smooth).  By taking a subsequence, we may assume that the $p_i$ approach a point $p\in S_\epsilon\cap X$ and that $T_{p_i}\big(G^{-1}G(p_i)\big)$ approaches a limit $\cal T$. If we let $M$ denote the stratum of our $a_G$ stratification which contains $p$, then we find that $T_pM\subseteq \cal T$ and ${\cal T}\subseteq T_pS_\epsilon$. As $S_\epsilon$ transversely intersects $M$, we have a contradiction, which proves the desired result.
\end{proof}

\begin{cor}\label{cor:main}
Suppose that $G$ is \L-analytic at $p$. Then, the Milnor fibrations for $G-G(p)$ inside a ball, centered at $p$, exist. 
\end{cor}
\begin{proof} As we saw in \secref{sec:milnorcondab}, this follows immediately from Milnor's condition's (a) and (b). 
\end{proof}

\smallskip

We include the following just to emphasize that \L-analytic have many properties that one associates with complex analytic functions.

\smallskip

\begin{cor} Suppose that $G$ is an \L-analytic function which is nowhere locally constant, and that $k\geq 2$. Then, $G$ is an open map.
\end{cor}
\begin{proof} This follows at once from \thmref{thm:main} and \corref{cor:open}. \end{proof}

\section{Comments and Questions}\label{sec:comments}

We have many basic questions.

\begin{ques} Are there real analytic maps which are \L-maps but do not have positive \L-weight?
\end{ques}

We suspect that the answer to the above question is ``yes''.

\medskip

Given \thmref{thm:lweightposcomp}, it is natural to ask:

\begin{ques} Is the composition of two \L-maps an \L-map? What if the maps are \L-analytic?
\end{ques}

\smallskip

Considering \exref{exm:mix}, we are led to ask:

\begin{ques} If one takes two simple \L-maps, in disjoint variables, into $\R^4$ or $\R^8$ and multiplies them using quaternionic or octonionic multiplication (analogous to the complex situation in \exref{exm:mix}),  does one obtain a new simple \L-map?
\end{ques}

We feel certain that the answer to the above question is ``yes'', but we have not actually verified this.

\bigskip

It is not difficult to produce examples of real analytic $f=(g, h)$ which possess Milnor fibrations inside balls, and yet are {\bf not} \L-analytic. It is our hope that the \L-analytic condition is strong enough to allow one to conclude that the Milnor fibration inside an open ball is equivalent to the Milnor fibration {\bf on the sphere} (the fibration  analogous to $f_{{}_\C}/|f_{{}_\C}|:S_\epsilon-S_\epsilon\cap X\rightarrow S^1\subseteq \C$, which exists when $f_{{}_\C}$ is complex analytic). 

As we discussed in the introduction, in the real analytic setting, $f/|f|$ does not necessarily yield the projection map of a fibration from $S_\epsilon-S_\epsilon\cap X$ to $S^1$; when this map is, in fact, a locally trivial fibration (for all sufficiently small $\epsilon$), one says that $f$ satisfies the {\bf strong Milnor condition} at the origin (see VII.2 of \cite{seadebook}).

It is not terribly difficult to try to mimic Milnor's proof that the Milnor fibration inside the ball and the one on the sphere are diffeomorphic; one wants to integrate an appropriate vector field. In the case where $f=(g, h)$ is real analytic and satisfies Milnor's conditions (a) and (b), one finds that one needs a further condition:

\begin{defn} Let $\omega:=-h\nabla g +g\nabla h$, and assume that $f(\0)=\0$. We say that {\bf $f$ satisfies Milnor's condition (c) at $\0$} if and only if there exists an open neighborhood $\W$ of $\0$ in $\U$ such that, at all points $\mathbf x$ of $\W-X$, if $\mathbf x$ is a linear combination of $\nabla g$ and $\nabla h$, then
$$
|\omega(\mathbf x)|^2\big(\mathbf x\cdot\nabla|f|^2(\mathbf x)\big)>\big(\nabla|f|^2(\mathbf x)\cdot\omega(\mathbf x)\big)\big(\mathbf x\cdot\omega(\mathbf x)\big).
$$
\end{defn}

\bigskip

It is not too difficult to show that, if $f$ is a simple \L-analytic function, there exists an open neighborhood $\W$ of $\0$ in $\U$ such that, at all points $\mathbf x$ of $\W-X$, $\mathbf x\cdot\nabla|f|^2(\mathbf x)>0$. It is also trivial that a simple \L-analytic function has $\nabla|f|^2(\mathbf x)\cdot\omega(\mathbf x)=0$. Milnor's condition (a) implies that $\W$ can be chosen so that $\omega(\mathbf x)\neq 0$ for $\mathbf x\in\W-X$. Thus, we conclude that simple \L-analytic functions satisfy the strong Milnor condition.

\smallskip

However, the real question is:

\begin{ques} Do general \L-analytic maps into $\R^2$ satisfy Milnor's condition (c), and so satisfy the strong Milnor condition? What about real analytic maps with positive \L-weight? In light of \propref{prop:jacequiv}, what about real analytic maps with positive \L-weight where $\nabla g$ and $\nabla h$ have comparable in magnitude?
\end{ques}

\medskip

Finally, having seen in this paper that \L-analytic pairs of functions share some important properties with complex analytic functions, one is led to ask a very general question:

\begin{ques} Do \L-analytic maps form an interesting class of functions to study for reasons having nothing to do with Milnor fibrations?
\end{ques}

\newpage
\bibliographystyle{plain}
\bibliography{Masseybib}
\end{document}